# SOME PROPERTIES OF TRANSFORMS IN CULTURAL THEORY

## BY PAUL BALLONOFF[1]

*Abstract:* It is shown that, in certain circumstances, systems of cultural rules may be represented by doubly stochastic matrices denoted $\Pi$, called "possibility transforms," and by certain real valued "possibility densities" $\pi=(\pi_1, \pi_2, ..., \pi_n)$ with inner product $\langle\pi,\pi\rangle=\Sigma_i\pi_i^2=1$. Using such objects we may characterize a certain problem of ethnographic and ethological description as a problem of prediction, in which observations are predicted by properties of fixed points of transforms of "pure systems", or by properties of convex combinations of such "pure systems". That is, ethnographic description is an application of the Birkhoff theorem regarding doubly stochastic matrices on a space whose vertices are permutations.

This paper follows from and adopts the definitions used in previous papers [4, 5] to which the reader is referred for details. In brief we assume a finite non-empty set $\mathbb{P}$ whose members are called <u>individuals</u>. Such $\mathbb{P}$ is organized into an <u>evolutionary structure</u> $\mathbb{S}$ as a quintuple ($\mathbb{P}$, $\mathbb{R}$, *D, B, M*) where $\mathbb{R}$ is a non-empty set of "rules", and *D, B,* and *M* are binary relations on $\mathbb{P}$. We assume that each evolutionary structure satisfies the following four axioms: (1) *D* is totally non-symmetric and transitive; (2) *M* is symmetric; (3) if b*D*c and there exists no d$\in\mathbb{P}$, d$\neq$b,c for which b*D*d and d*D*c, then we write c*P*b, and then require B = { (b,c) | there exists d$\in\mathbb{P}$ with both d*P*b and d*P*c }; and (4) #b$M \leq 2$, where b$M$:={c$\in\mathbb{P}$ | (b,c) $\in M$}. A <u>rule</u> R$\in\mathbb{R}$ is a statement concerning the relationships between the *D, B,* and/or *M*, which does not violate those four axioms. Given a family of subsets $\mathbb{G}^t\in\mathbb{P}$, indexed by a set t$\in\mathcal{T}$ of consecutive non-negative integers starting with 0, each $\mathbb{G}^t$ is thought of as the <u>generation</u> at time t. Then $\mathcal{G}=\{\mathbb{G}^t \mid t\in\mathcal{T}\}$ is called a <u>descent sequence</u> of $\mathbb{S}$ in case, for all $\mathbb{G}^t\in\mathcal{G}$, each cell B occurs in only one generation, each subset M occurs in only one generation, and when $\mathbb{G}^t\in\mathcal{G}$, b$\in\mathbb{G}^t$, and c*P*b, then c$\in\mathbb{G}^{t-1}$ (that is, the set $\mathbb{G}^{t+1}$ contains all of, and only, the immediate descendants of individuals in $\mathbb{G}^t$). Thus the generations of each descent sequence $\mathcal{G}$ partitions its total population. We assume a "Darwinian Sequences axiom" which says that all descent sequences of a given evolutionary structure can be traced back through a chain of descent in an unbroken series of non-empty generations, to the same date of initial origin.

We define by enumeration a set of closed (cyclic) objects which we call <u>regular structures</u>. We illustrate using a dot to represent an individual, a circle around two dots (say, b and c) shows b*M*c, a line between two dots (say, d,e) shows d*B*e. We give the simple closed cycles names Mn, where n = the number of M-sets in each figure. For example a diagram of relations and resulting M2 regular structure is:

[1] I thank Dick Greechie for his generous comments on earlier versions of this text. All errors are my own.

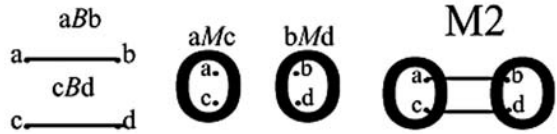

and

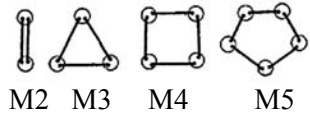

M2  M3  M4  M5

are the regular structures with 2, 3, 4 and 5 *M*-sets. All individuals in a particular regular structure are necessarily of the same generation. Each generation is composed of a set of regular structures. An ordered list counting the numbers of regular structures present in a particular $\mathbb{G}^t$ is a <u>configuration</u> C:=$(m_1, ..., m_j, ...)$ where the coefficient $m_j$ is the number of regular structures of type Mj in C. Thus a configuration consisting only of 2 of the M2 structures would be written (0,0,2,0, …). We shall denote that $\mathbb{C} = \{ C_i \mid i = 1,2, ..., n\}$ is a finite non-empty set of n distinct configurations $C_i$. We here consider only finite sets $\mathbb{C}$. In general, if $C_i$ and $C_j$ are configurations then $C_i+C_j$ is also a configuration, though $C_i+C_j \in \mathbb{C}$ is not required (since $\mathbb{C}$ is finite). If C is a configuration, then $\mu(C)=\Sigma(jm_j)$, or simply "$\mu$" when C is understood, is the number of reproducing marriages in C (the "marriage number of $C_i$"); $\beta$ is the number of cells induced by *B* (sibships) in C; and $\gamma=2\mu$ is the total population of the generation $\mathbb{G}^t$ on which $C_i$ is formed. If $\mathcal{G}$ is a descent sequence and, for $t \in \mathcal{T}$, $\mathbb{G}^t \in \mathcal{G}$ is a generation in $\mathcal{G}$ of size $\gamma^t$, $\mu^t$ is the number of marriages on $\mathbb{G}^t$, and $\beta^t$ the number of sibships on $\mathbb{G}^t$, then (when defined) $\mu^{t-1}=\beta^t$ and $\mu^{t+1} \leq \mu^t$, since a generation can not have more sibships (which arise by descent) than are reproducing marriages in its predecessor generation.

We define $\mathcal{P}(\mathbb{C}):=\{\xi \mid \xi \subseteq \mathbb{C}\}$. Let $\xi \in \mathcal{P}(\mathbb{C})$. We let each $\xi \in \mathcal{P}(\mathbb{C})$ correspond to a unique "row vector" $(\xi_1, \xi_2, ..., \xi_n)$ where $\xi_j=1$ if $C_i \in \xi$, and $\xi_i=0$ otherwise, called the <u>content list</u> of $\xi$. We denote $\Xi:=\{(\xi_1, \xi_2, ..., \xi_n) \mid \xi \in \mathcal{P}(\mathbb{C}), \xi=(\xi_1, \xi_2, ..., \xi_n)$ is the content list of $\xi\}$. For simplicity by abuse of notation we use the notation $\xi \in \mathcal{P}(\mathbb{C})$ and also $\xi \in \Xi$ for the corresponding content list $(\xi_1, \xi_2, ..., \xi_n)$. We denote the empty configuration as $C_0$ but list only non-empty configurations. For example, if n=3 so that $\mathbb{C}=\{\xi_1, \xi_2, \xi_3\}$, and if $\xi=\{\xi_1, \xi_2\}$ then the content list $\xi= (1,1,0)$. The transpose of a vector $\xi=(\xi_1, \xi_2, ..., \xi_n)$ is denoted $\xi^T=(\xi_1, \xi_2, ..., \xi_n)^T$. The objects $\xi$ and $\xi^T$ represent the same set of configurations. Let $C_i \in \mathbb{C}$ and let $\mu_C=s$, s an integer >0, be the marriage number of $C_i$. Then: $\mathbb{C}_s:=\{ C_i \mid C_i \in \mathbb{C}$ and $\mu_C=s\}$ is a <u>set of configurations of order s</u>, and $\mathcal{P}_s:=\mathcal{P}(\mathbb{C}_s)$ is a set of subsets of $\mathbb{C}_s$.

A rule can be represented by defining a transform $\mathbf{R}=[r_{ij}]$ where $r_{ij}=1$ if the rule allows a transition from $C_j$ to $C_j$ otherwise $r_{ij}=0$. We arrange the $r_{ij}$ on the ordered pairs $(C_i, C_j)$, $1 \leq i,j \leq n$, in rows $r_i=(r_{i1}, r_{i2},...,r_{in})$, $0 \leq i,j \leq n$ in the standard order; there are n such rows.

If $\mathbf{R}=[r_{ij}]$ is a transform and $r_i=(r_{i1}, r_{i2}, \ldots, r_{in})$ is a row of $\mathbf{R}$, then $r_i$ shows which transitions, from each $\xi_j$ to $\xi_i$ are allowed by $\mathbf{R}$.

> *Definition 1*: Let $\mathbb{S}$ be an evolutionary structure with a finite non-empty set $\mathbb{C}$ of configurations on generations of $\mathbb{S}$, and let $\check{\mathbf{R}}$ be a finite non-empty set of transforms on $\mathbb{C}$. Then
> 1. Let $\hat{\mathbf{H}}$ be the free semigroup generated by $\check{\mathbf{R}}$, called the set of <u>transforms generated by $\check{\mathbf{R}}$</u>;
> 2. If $\alpha \in \hat{\mathbf{H}}$ then $\alpha$ is called a <u>history generated by $\check{\mathbf{R}}$</u>

A <u>history</u> is thus a string of one or more transforms applied in succession on a descent sequence. For example, if $\check{\mathbf{R}}$ is a set of transforms and $\mathbf{Q},\mathbf{R},\mathbf{S} \in \check{\mathbf{R}}$ then $\alpha=\mathbf{QRS}$ is a string of transforms of $\check{\mathbf{R}}$ and thus $\alpha$ is a history. Histories are associative under composition by concatenation. We do not allow that an empty chain is a history. Histories occur on an evolutionary structure $\mathbb{S}$, so existence of a history (that is, of application of certain rules in a stated sequence) implies existence of some non-empty descent sequence $\mathcal{G}$ of $\mathbb{S}$ on whose generations the rules are applied in that order. If $\hat{\mathbf{H}}$ is the set of histories generated by a set $\check{\mathbf{R}}$ of transforms then $\check{\mathbf{R}} \subseteq \hat{\mathbf{H}}$ and if $\mathbf{R} \in \check{\mathbf{R}}$, $\mathbf{R} \in \hat{\mathbf{H}}$.

*Notational Conventions:* Hereafter we adopt the conventions that $\mathbb{C}$ is a finite non-empty set of $n \geq 1$ non-empty configurations, $\mathcal{P}(\mathbb{C})$ is a set of all subsets of $\mathbb{C}$, $\Xi$ the set of contents lists for $\xi \in \mathcal{P}(\mathbb{C})$, $\check{\mathbf{R}}$ a non-empty set of transforms, $\hat{\mathbf{H}}$ a set of transforms generated by $\check{\mathbf{R}}$, and we let $\alpha \in \hat{\mathbf{H}}$, $\alpha=[\alpha_{ij}]$ $0 \leq i \leq n$, $0 \leq j \leq n$ and for $\xi,\phi \in \Xi$, $\alpha\xi^T=\phi^T$.

> *Definition 2*: Let $\mathbb{C}$ be a non-empty set of configurations, let $\check{\mathbf{R}}$ be a non-empty set of transforms on $\mathbb{C}$, let $\mathbf{R},\mathbf{S} \in \check{\mathbf{R}}$, and let $\mathbf{R}=[r_{ij}]$ and $\mathbf{S}=[s_{ij}]$. Let $\alpha=\mathbf{SR}$. Let $r_i=(r_{i1}, \ldots, r_{in})$ be the $i^{th}$ row of $\mathbf{R}$ and $s_j=(s_{1j}, \ldots, s_{nj})$ be the $j^{th}$ column of $\mathbf{S}$, where $0 \leq i,j \leq n$. Then
> $$\alpha_{ij}=r_i s_j=(r_{i1}s_{1j} \circ r_{i2}s_{2j} \circ \ldots \circ r_{iv}s_{vj})$$
> where $xy = \min\{x,y\}$, and $x \circ y = \max\{x,y\}$.

Then $\alpha=[\alpha_{ij}]$ is the transform of $\alpha$, which we call the <u>logical product</u> (by matrix multiplication) of the transforms $\mathbf{S}$ and $\mathbf{R}$. Definition 2 simply incorporates the result of Theorem 3 of [1]. We apply the same logical product multiplication to the product of transforms with vectors $\xi \in \Xi$. Thus if $\mathbb{C}$ is a set of configurations, $\xi \in \Xi$ then $\alpha\xi^T=\phi^T$, where also $\phi \in \Xi$.

For example, let $\mathbb{C}$ be a set of configurations, $\alpha=\mathbf{RS}$ be a history acting on configurations of $\mathbb{C}$, and $\alpha=[\alpha_{ij}]$. Then $\alpha_{ij}=0$ if $\mathbf{S}$ allows no transition to $C_i$ from any $C_j \in \mathbb{C}$; while $\alpha_{ij}=1$ if $\mathbf{S}$ allows a transition from any $C_j \in \mathbb{C}$ to $C_i \in \mathbb{C}$, and $\mathbf{R}$ allows a transition to $C_j$ from at least one of the $C_k$ allowed under $\mathbf{S}$. This operation is associative, so given a history (a string of transforms) $\alpha$, using Definition 2 we can construct a single transform that represents the effect of $\alpha$.

*Definition 3*: Under the notational conventions: (a) A history $\alpha$ is <u>viable</u> if there exists $\xi \in \Xi$, $\xi > 0$, $\Sigma_i \xi_i \geq 1$ such that $\alpha \xi^T = \xi^T$, and for $\varphi \in \Xi$, $\varphi > 0$ and $\varphi \leq \xi$ then $\alpha \varphi^T = \varphi^T$. Then we also say that $\alpha$ is <u>viable on</u> $\xi$. (b) A <u>minimal structure</u> of $\alpha$ is a non-empty configuration C such that $\xi_C = 1$ in $\xi$, $\alpha$ is viable on $\xi$, and if there exists a nonempty configuration D such that when $\varphi \in \Xi$, $\varphi > 0$, $\alpha$ is viable on $\varphi$ and $\varphi_D = 1$ then $\Sigma(jc_j) \leq \Sigma(jd_j)$.

If $\alpha$ is viable on $\xi$ then $\alpha$ is viable on any $\varphi \in \Xi$, $\varphi \leq \xi$, provided $\varphi > 0$. If $\alpha$ is viable, then $\alpha$ has at least one minimal structure. If C is a minimal structure of $\alpha$, then $\Sigma(jc_j) = s$ is the <u>structural number</u> of $\alpha$. Clearly, every minimal structure of a history $\alpha$ has the same structural number, and every such minimal structure is a "fixed point" of $\alpha$ (see also [4]). A history $\alpha$ acting only on a minimal structure of $\alpha$ is viable.

Let $\mathbb{C}$ be a non-empty set of configurations, and let $\mathcal{T}$ be the set of all transforms on $\mathbb{C}$ which are allowed by the rules of construction of configurations. Then we say that $\mathcal{T}$ is a <u>full</u> set of transforms on $\mathbb{C}$. Let $\alpha = [\alpha_{ij}]$ be a transform and let $\alpha^T = [\beta_{ji}]$ where $\beta_{ji} = \alpha_{ij}$ be a matrix which is the transpose of $\alpha$. Such $\alpha^T$ is also seemingly a transform, provided there is a rule generating $\alpha^T$, which "reverses" the action of $\alpha$, since $\beta_{ji} = \alpha_{ij}$; the object $\alpha^T$ would thus appear to allow j to transform to i, if $\alpha$ allows i to transform to j. However, if we again let $\mu_C = \Sigma(jc_j)$, $\mu_D = \Sigma(jd_j)$ and $\mu_D > \mu_C$, then $\alpha_{DC} = 0$ even if $\alpha_{CD} = 1$, because under the rules of construction of the *M* and *B* relations, there can not be more cells defined by the relation *B* created in one generation, than is the marriage number of the preceding generation. Thus, there can be no rule creating such $\alpha^T$ and thus not every matrix which is a transpose of a transform, is itself a transform allowed by the rules of construction of configurations. For similar reasons, the inverse of a transform in a full set is not a necessary member of that set. Thus a full set of transforms in general is also not a group.

Note that each $\xi_i$ (respectively $\phi_i$) is a number either 0 or 1, so if $\alpha \xi^T = \phi^T$ and $\xi, \phi \neq 0$, then $\Sigma \xi_i \geq 1$ and $\Sigma \phi_i \geq 1$. If $\Sigma \phi_i = 1$ then $\alpha$ has exactly one possible outcome, and if $\Sigma \xi_i = 1$ then $\alpha$ acted on a set consisting of only one initial configuration. We wish to describe the "relative possibilities" when more than one possible outcome might occur. This problem includes the following: how to describe the "relative possibilities" of outcomes when more than one immediately previous configuration might have existed, and we do not know which intermediate configuration actually occurred.

*Definition 4*
1. If $\alpha \in \hat{\mathbf{H}}$ is a history, a <u>possibility transform of</u> $\alpha$ is $\Pi = [p_{ij}]$ in which $0 \leq p_{ij} \leq 1$, $\Sigma_j p_{ij} \leq 1$, $p_{ij} > 0$ iff $\alpha_{ij} = 1$. Then $\mathbf{\Pi} := \{ \Pi \mid \Pi \text{ is a possibility transform}\}$ is a set of possibility transforms.
2. Let $\mathbf{\Pi}$ be a non-empty set of possibility transforms of dimension $n > 0$, let $\Pi, \Theta \in \mathbf{\Pi}$, and let $\xi, \omega \in \Xi$, where $\Sigma_i \xi_i = w \leq n$ (respectively $\Sigma_i \omega_i = x \leq n$). Let the products $\Pi \Theta$ be computed using ordinary arithmetic product of two matrices, and the products $\xi \Pi$ (respectively $\Pi \xi^T$) as the arithmetic product of a vector and a matrix. Noticing that $\Sigma_i \xi_i = w$, $w \leq n$, then the <u>possibility density</u> $\pi$ of $\xi \Pi$ (respectively $\Pi\omega^T$) is $\pi:=(\pi_1, \pi_2, ...\pi_n)$ (respectively $\omega^T:=(\omega_1, \omega_2, ...\omega_n)^T$) where $\pi_i=\Sigma_j(p_{ij}\xi_j/w)$ (respectively $\omega_i=\Sigma_j(p_{ji}\omega_j/w)$.
3. If $\pi,\omega$ are possibility densities then $<\pi,\omega>=\Sigma_i(\pi_i\omega_i)$ is the <u>inner product</u> of $\pi$ and $\omega$.

*Comment 1*: Since each $\xi_j=0$ or $\xi_j=1$, and each $p_{ij}\xi_j>0$ iff $\xi_j=1$, then if $\xi_j=1$, $p_{ij}\xi_j=p_{ij}$, otherwise $p_{ij}=0$. So while Definition 4.1 allows $\Sigma_j(p_{ij}\xi_j)\leq 1$, $\Sigma_j(p_{ij}\xi_j)=1$ iff $\Sigma_j p_{ij}=1$, which implies $\Sigma\xi_j>0$.

A possibility density $\pi$ of $\xi\Pi$ answers "if we know the transition ends in one of the non-empty configurations of $\xi$ and the descent sequence applies the possibility transform $\Pi$, what if any are the non-empty configurations from which it may have started, and the relative possibility of each?" A possibility density $\omega^T$ of $\Pi\xi^T$ answers "if we know the transition starts in one of the non-empty configurations of $\xi$ and the descent sequence applied the possibility transform $\Pi$, what are the non-empty configurations from which the transition might end and the relative possibility of each?" The inner product answers both questions at once. Given these interpretations, it is reasonable to impose this axiom:

*Axiom 1:* If $\pi=(\pi_1, \pi_2, ...\pi_n)$ is a possibility density then $\Sigma_i\pi_i\leq 1$.

Since a viable history must have a non-empty descent sequence, we would thus like to know when $<\pi,\omega>=1$.

*Theorem 1:* Under conditions of Definitions 3 and 4 and the notational conventions, if $\alpha,\beta\in\hat{\mathbf{H}}$, $\Pi,\Theta\in\Pi$ are possibility transforms of $\alpha$, $\beta$ respectively, $\xi,\phi\in\Xi$, $\pi$ the possibility density of $\xi\Pi$, $\omega^T$ the possibility density of $\Theta\phi^T$, $\Sigma_i\xi_i=w$, all these lists and matrices have dimension n, then $<\pi,\omega>=1$ iff all of:
(i) $\phi\neq 0$; (ii) $\xi\neq 0$; (iii) for each i, $\Sigma_j p_{ij}=1$ and $\Sigma_j q_{ij}=1$; (iv) $\phi=\xi$; (v) $\Sigma\xi_j=\Sigma\phi_j=w$.

*Proof:* The products $\xi\Pi$ and $\Theta\phi^T$ are defined since the matrices and vectors all have the same dimension, and thus the respective possibility densities $\pi$, $\omega^T$ are defined. Note that $<\pi,\omega^T>=\Sigma_i((\Sigma_j(p_{ij}\xi_j/w)(\Sigma_j(q_{ij}\phi_j/w))$. Conditions (i) and (ii) are required since if either $\phi,\xi=0$, $<\pi,\omega^T>=0$, since all of the products then involve at least one term $= 0$, thus all products $= 0$. Condition (iv) results since if there is any i such that $\phi_i=1$ but $\xi_i=0$, or any $\phi_i=0$ but $\xi_i=1$, then there will be a product $p_{ij}\xi_i q_{ji}\phi_i=0$ in a numerator, but for which not both $p_{ij}=0$ and $q_{ij}=0$. But in that case even if $\Sigma_j p_{ij}=1$ and $\Sigma_j q_{ij}=1$, at least one of the $p_{ij}\neq 0$ or $q_{ij}\neq 0$ values will be disregarded in computing sums, so $\Sigma_i((\Sigma_j(p_{ij}\xi_j/w)(\Sigma_j(q_{ij}\phi_j/w))<1$. Thus necessarily $\phi=\xi$ so (iv) is shown. To show conditions (iii) recall that $0\leq p_{ij}\leq 1$ and $0\leq q_{ij}\leq 1$, $\Sigma_j p_{ij}\leq 1$, $\Sigma_j q_{ij}\leq 1$, and thus every sum $\Sigma_j(p_{ij}\xi_i/w)\leq 1$, and every $\Sigma_j(q_{ij}\phi_i/w)\leq 1$. And if not both: (1) $\Sigma_j p_{ij}=1$ for all i then there is at least one i such that $\Sigma_j(p_{ij}\xi_{ij}/w)<1$, and (2) $\Sigma_j q_{ij}=1$ for all i then there is at least one i such that $\Sigma_j(q_{ij}\phi_i/w)<1$; so $<\pi,\omega^T>=1$ only if both $\Sigma_j p_{ij}=1$ and $\Sigma_j q_{ij}=1$. Thus (iv) is met, and since $\phi=\xi$ then $\Sigma\xi_j=\Sigma\phi_j=w$, so (v) is met. Thus only if (i) through (v) then $<\pi,\omega>=1$. But also by similar arguments, if $<\pi,\omega>=1$ then conditions (i) through (v) are met.//

If $\alpha$ is viable then conditions (i), (ii), (iv) and (v) of Theorem 1 are met, but these conditions are requisite for (iii), so a history $\alpha$ is viable if its transform $\Pi$ meets Theorem 1. From Theorem 1 (iv) $\phi=\xi$ and Comment 1, if any row or column of a possibility transform $\Pi$ or $\Theta$ is that of some $\phi_i=\xi_i=0$ then the corresponding $p_{ij}=0$ and $q_{ij}=0$. Thus when discussing possibility transforms meeting Theorem 1, we can consider such possibility transform to be only its non-zero rows and columns, and the corresponding lists only the non-zero entries, and all thus of dimension w; we call such forms <u>reduced</u>.

*Lemma 1:* If $\pi=(\pi_1, \pi_2, ..., \pi_n)$ and $\omega^T=(\omega_1, \omega_2, ..., \omega_n)^T$ are possibility densities meeting Theorem 1, then $\Sigma_i\pi_i=1$ and $\Sigma_i\omega_i=1$.

*Proof:* From Theorem 1(iii) each non zero row of $\Pi$ is such that $\Sigma_j p_{ij}=1$, $0\leq p_{ij}\leq 0$. From Definition 4.2, using the multiplication defined in 4.3, and from Theorem 1 (v), each $\omega$ is a convex combination of w rows of $\Pi$, each of which row is weighted by 1/w. Thus, $\Sigma_i\omega_i=w(1/w)1=1$. By a similar argument, each $\Sigma_i\xi_i$ is a sum of exactly all non-zero entries $p_{ij}$, with exactly the same weights since there are w non-zero columns, and the denominators are all w, and, though possibly taken in a different sequence, that sum is thus given by $\Sigma_i\pi_i=\Sigma_i\omega_i=1$.//

A <u>doubly stochastic matrix</u> is a matrix, each of whose rows and columns are non-negative numbers that sum to 1.

*Theorem 2:* A possibility transform $\Pi$ meeting Lemma 1 and Theorem 1 is doubly stochastic on those rows and columns for which $\phi_i=\xi_i=1$

*Proof:* Let $\Pi$ be a possibility transform satisfying Theorem 1. We consider only rows and columns for which $\phi_i=\xi_i=1$ (we need consider only the reduced form of $\Pi$). There are w such rows or columns. Theorem 1 establishes the result as to rows of $\Pi$. We extend the same reasoning as in Lemma 1. Given a non-empty set of configurations $\xi$, the possibility density $\omega$ is a sum of w columns, each weighted by 1/w, whose sum is w, so the sum of each such column is $w(1/w)=1$.//

*Theorem 3:* Let $\alpha$ be a history, let $\xi\in\Xi$, let $\alpha$ be viable on $\xi$, let be $\Pi$ the possibility transform of $\alpha$, let and $\pi$ be the possibility density of $\xi\Pi$ then $<\pi,\pi^T>=\Sigma_i\pi_i^2=1$.

*Proof:* Obvious. We note Theorem 1, set $\Theta=\Pi$, set $\phi=\xi$, then find the possibility densities $\pi$ of $\xi\Pi$ and $\omega^T$ of $\Pi\xi^T$. But since we assume $\Pi$ is symmetric, then $\pi^T=\omega^T$ so substituting in Theorem 1 gives $<\pi,\pi>=1$. Since $\pi=(\pi_1, \pi_2, ..., \pi_n)$, then also $<\pi,\pi>=\Sigma_i\pi_i^2=1$.//

An important example of a set $\xi$ of configurations allowing to meet Theorems 1 through 3 is a set of the minimal structures of a transform (history, rule) $\alpha$.

*Definition 5*: Let $\mathbb{C}_s$ be a non-empty set of configurations of order s, and let $\mathscr{P}_s(\mathbb{C}_s)=\{\xi \mid \xi \subseteq \mathbb{C}_s\}$. Let $\hat{\mathbf{H}}_s$ be the full set of transforms on $\mathbb{C}_s$. Let $\Pi_s:=\{\Pi_\alpha \mid \alpha \in \hat{\mathbf{H}}_s, \Pi_\alpha$ is the possibility transform for $\alpha\}$.
1. Then $S_s:=(\mathbb{C}_s, \mathscr{P}_s(\mathbb{C}_s), \hat{\mathbf{H}}_s, \Pi_s)$ is a <u>cultural structure of order s</u>.
2. If $S_s$ is a cultural structure of order s, $\alpha \in \hat{\mathbf{H}}_s$ a history with minimal structure $\xi_m \in \mathbb{C}_s$ and structural number s, $\alpha=[\alpha_{ij}]$, $\alpha_{ij}=1$ iff i=j=m, and $\hat{\mathbf{H}}_s=\{\alpha\}$, then $S_s$ is a <u>pure system</u> of $\alpha$.

A pure system of $\alpha$ is viable. Let $\alpha$ be a history with structural number s>0, and $S_s$ a pure system of $\alpha$, then $\hat{\mathbf{H}}_s=\{\alpha\}$, and since $\alpha$ acts only on a particular minimal structure of $\alpha$, then in a pure system $\alpha=\alpha^{-1}=\alpha^2$. In a pure system, the possibility transform $\Pi$ of $\alpha$ has but one entry $p_{ii}=1$, where i indexes the minimal structure configuration $\xi_m$ and all other $p_{jk}=0$, jk≠ii. Then for such $\Pi$, tr$\Pi$=1, indeed tr$\Pi\Pi$=1, and such $\Pi$ is obviously symmetric.

*Definition 6:* Under conditions of Definition 4 let $\hat{\mathbf{H}}$ be a set of transforms, let $\alpha \in \hat{\mathbf{H}}$, let $\Pi_\alpha$ be the possibility transform of $\alpha$, and let
$$\Pi = \{\Pi_\alpha \mid \alpha \in \hat{\mathbf{H}}, \text{ and } \Pi_\alpha \text{ is the possibility transform of } \alpha\}.$$
Let $v_\alpha$ be a real number $0 \leq v_\alpha \leq 1$ such that $\Sigma_\alpha v_\alpha = 1$. Then $\Psi := \Sigma_\alpha v_\alpha \Pi_\alpha$, is a <u>convex combination</u> of possibility transforms

Let $\theta \subseteq \hat{\mathbf{H}}$ be a non-empty subset of $\hat{\mathbf{H}}$ such that $v_\alpha > 0$ iff $\alpha \in \theta$. Now let $\Psi=[y_{ij}]$. Recall from Definition 4.1 that if $\alpha$ is a history, the possibility transform of $\alpha$ is $\Pi_\alpha=[p_{ij}]$ in which (1) $0 \leq p_{ij} \leq 1$, (2) $\Sigma_j p_{ij} \leq 1$, and (3) $p_{ij} > 0$ iff $\alpha_{ij}=1$. Clearly such $\Psi$ meets (1) and (2). The equalities occur in (2) only if $\alpha$ is viable. If we assume the set $\Pi$ of Definition 6 includes only possibility transforms of pure systems, then clearly also, tr$\Psi$=1, and also tr$(\Sigma_\alpha v_\alpha \Pi_\alpha)$=tr$\Psi$=1. Then for each such pure system, if $\Pi_\alpha=[p_{ij}] \in \theta$ there exists an $\alpha$ such that exactly one $\alpha_{ij}=1$. And therefore if $y_{ij}>0$ there also exists at least one $\alpha$ with at least one $\alpha_{ij}=1$. This shows that when $\Psi = \Sigma_\alpha v_\alpha \Pi_\alpha$, $\alpha \in \hat{\mathbf{H}}$, $\Pi_\alpha \in \Pi$, then also $\Psi \in \Pi$. We may then need to construct a corresponding history $\alpha$ (which is what anthropologists typically do in creating a minimally structured genealogy to describe the action of a "cultural rule"), but we know that (at least) one exists.

Consider an ethnologist (or field biologist) observing a species of individuals in a descent sequence of some evolutionary structure. At the outset of work, the ethnologist hypothesizes that tr$\Psi$=1 meaning, that the descent sequence has some marriage rule(s), which can be characterized by some pure system(s). After observing, the ethnologist claims a subset of rules $\theta \subseteq \hat{\mathbf{H}}_s$ such that tr$\Psi$=tr$(\Sigma_\alpha v_\alpha \Pi_\alpha)$=1 for $\alpha \in \theta$. Of course $\theta$ may consist of a single rule $\alpha$ (for which therefore $v_\alpha$=1). It is known [1, 2, 3], that that if a descent sequence follows a rule $\alpha$ with structural number s, that certain "demographic" measures can be computed from the convex combination of structural numbers of the rules thus observed. Thus the ethnologist can perform several simultaneous empirical observations on the descent sequence at date t, including: ask if they observed at least one rule $\alpha$ with structural number s>0 (that is, observed a non-empty set $\theta \subseteq \hat{H}$ such that tr$\Psi$=tr($\Sigma_\alpha v_\alpha \Pi_\alpha$)=1, for $\alpha \in \theta$); and ask if the statistics are those predicted from the structural number or convex combination of structural numbers of the claimed rule(s). Where other features such as "kinship terminologies" are also studied, the ethnologist can also ask if the kinship terminology of each $\alpha \in \theta \subseteq \hat{H}$ maps "properly" to a descent sequence composed of a pure system using the minimal structure of $\alpha$. Ethnographic description thus amounts to asking which if any transform(s) $\alpha$ have a fixed point that describes the empirically observable properties of a descent sequence following a rule.

The fact that the possibility densities of viable histories are doubly stochastic matrices opens to culture theory the application of the "Birkhoff theorem" [6], that the set of doubly stochastic matrices of order *n* is the convex closure of the set of permutation matrices of the same order, and the vertices (extreme points) of that set are those permutation matrices. At least many marriage and kinship systems are known to be described as permutations [7, 8]. Since pure systems are a very simple example of permutations, we have used only a very basic form of that result. The breadth of application this suggests for culture theory is large, and we shall exploit it in later papers.